\documentclass{article}
\usepackage{amsmath, amstext, amssymb}
\begin{document}
\centerline{\bf ON PRIMITIVE WEIRD NUMBERS OF THE FORM $\boldsymbol{2^kpq}$.}
\vskip 12pt\centerline{\it Douglas E. Iannucci}
\vskip 24pt\centerline{\parbox[t]{280 pt}{{\bf Abstract}\ \ We say a natural number~$n$ is abundant if $\sigma(n)>2n$, where $\sigma(n)$ denotes the sum of the divisors of~$n$. The aliquot parts of~$n$ are those divisors less than~$n$, and we say that an abundant number~$n$ is pseudoperfect if there is some subset of the aliquot parts of~$n$ which sum to~$n$. We say~$n$ is weird if~$n$ is abundant but not pseudoperfect. We call a weird number~$n$ primitive if none of its aliquot parts are weird. We find all primitive weird numbers of the form $2^kpq$ ($p<q$ being odd primes) for $1\le k\le14$. We also find primitive weird numbers of the same form, larger than any previously published. }}
\vskip 24pt\noindent{\bf 1. Introduction}
\vskip 12pt\noindent
We define  the arithmetic function $\sigma:\mathbb{N}\rightarrow\mathbb{N}$ by
$$\sigma(n)=\sum_{d\mid n}d,$$
the sum taken over the positive divisors~$d$ of~$n$. We say~$n$ is {\it abundant\/} if $\sigma(n)>2n$, {\it perfect\/} if $\sigma(n)=2n$, and {\it deficient\/} if $\sigma(n)<2n$. 

The {\it aliquot parts\/} of~$n$ are the positive divisors of~$n$ that are less than~$n$. We may denote this set by 
$$\mathcal{A}_n=\{d\mid n:0<d<n\}.$$
We say that~$n$ is {\it pseudoperfect\/} if~$n$ is abundant and there exists a subset $S\subset\mathcal{A}_n$ such that
$$n=\sum_{d\in S}d.$$
We say~$n$ is {\it weird\/} if~$n$ is abundant but not pseudoperfect. The smallest weird number is~70, and all known weird numbers are even. Benkoski and Erd\H{o}s~[1] have shown that the weird numbers have positive asymptotic density, {\it a fortiori\/} there exist infinitely many of them. However, their infinitude also follows from the simple fact that $pw$ is weird if~$w$ is weird and $p>\sigma(w)$ is a prime. 

We say~$n$ is {\it primitive weird\/} if~$n$ is weird and none of its aliquot parts are weird. It is unknown if infinitely many primitive weird numbers exist. 

Clearly all prime powers are deficient, and it is not hard to show that $2^kp^m$ is either deficient, perfect, or pseudoperfect for all $k\ge1$, $m\ge1$, and odd primes~$p$. Thus the simplest form a weird number may assume is $2^kpq$ for some~$k\ge1$ and odd primes $p<q$. Such numbers must be primitive weird. Pajunen~[4] listed all weird $2^kpq$ up to $k=8$. Kravitz~[3] searched for and discovered large numbers of this type, the largest being
$$ 2^{56}\cdot153722867280912929\cdot2305843009213693951.$$ 
We expand Pajunen's list to $k=14$, and we list several weird numbers $2^kpq$ of over a thousand digits' length. 

\vskip 24pt\noindent{\bf 2. Preliminaries}
\vskip 12pt\noindent
Consider a finite set $S\subset\mathbb{N}$ given by
$$S=\{a_1,\,a_2,\,a_3,\,\dots,\,a_n\}.$$
We define $\Sigma(S)$ by
$$\Sigma(S)=\left\{ \textstyle{ \sum_{j=1}^n\epsilon_ja_j\,:\,\epsilon_j=0\;\text{or}\;1,\;1\le j\le n}\right\}.$$
That is, $\Sigma(S)$ is the set of all possible sums taken from the elements of $S$, if no element is taken more than once. Clearly $|\Sigma(S)|\le2^n$. 

Thus~$n$ is pseudoperfect if and only if~ $n$ is abundant and $n\in\Sigma(\mathcal{A}_n)$. We define
$$a(n)=\sigma(n)-2n,$$
hence~$n$ is abundant if and only if $a(n)>0$. Since 
$$a(n)+n=\sigma(n)-n=\sum_{d\in\mathcal{A}_n}d,$$
it follows immediately that an abundant number~$n$ is pseudoperfect if and only if $a(n)\in\Sigma(\mathcal{A}_n)$. 
Fix~$k\ge1$ and let $n=2^kpq$, where $p<q$ are odd primes. Let $M=2^{k+1}-1$ and let $a=a(n)$, so that
$$a=M(p+1)(q+1)-(M+1)pq=M+Mp+Mq-pq.$$
Assuming~$n$ is abundant ({\it i.e.,\/} $a>0$), we have
\begin{equation}
a+pq=M+Mp+Mq\,;
\end{equation}
this is equivalent to
\begin{equation}
(p-M)(q-M)=M(M+1)-a.
\end{equation}
\vfill\eject
%
%
%
%
\noindent{\bf 3. Necessary and sufficient conditions for weird $\boldsymbol{2^kpq}$}
\vskip 12pt\noindent
Suppose $n=2^kpq$ is weird. We can show
\begin{equation}
a>M,
\end{equation}
\begin{equation}
M<p<2M.
\end{equation} 
For, if $a\le M$ we would have
$$a=\sum_{j\in J}2^j$$
for some $J\subset\{0,1,\dots,k\}$, whence $a\in\Sigma(\mathcal{A}_n)$. Similarly if $p\le M$ we would have
$$n=(2^k-1)pq+pq=\sum_{j=0}^{k-1}2^jpq+\sum_{j\in J}2^jq$$
for some $J\subset\{0,1,\dots,k\}$. Since we assume $p<q$, the left-hand side of~(2) is positive and
$$(p-M)<\sqrt{M(M+1)-a}<M,$$
whence $p<2M$. 
By~(2) and~(3) we have $4\mid a$ and 
\begin{equation}
M<a<M(M+1).
\end{equation}
Since~$n$ is weird, by~(1) we have $pq\ne r+sp+tq$ for all positive integers $r$, $s$, and $t\le M$. 

Conversely, fixing $k\ge1$, setting $M=2^{k+1}-1$, and given~$a$ as in~(5), we factor
$$de=\frac{M(M+1)-a}4,\qquad d<e.$$
If $p=M+2d$ and $q=M+2e$ are both prime, and if $pq\ne r+sp+tq$ for all positive integers $r$, $s$, and $t\le M$, then it follows from~(1) that $a\notin\Sigma(\mathcal{A}_{2^kpq})$, and hence $2^kpq$ is weird. 

Thus, for each fixed~$k\ge1$, we are able to conduct a finite search for all weird numbers $2^kpq$. We list all such for $1\le k\le14$, in Tables~1 to~41, which are found following the references. The number of weird numbers $2^kpq$, for each~$k$, $1\le k\le14$, are as follows.
\vskip 8pt
\begin{tabular}{llllllll}
$k$ & \# of weirds & $k$ & \# of weirds & $k$ & \# of weirds & $k$ & \# of weirds \\
1 & 1 & 5 & 10 & 9&115& 12& 683\\
2 & 1 & 6 & 23 & 10 & 210 & 13 & 1389\\
3 & 5 & 7 & 29 & 11 & 394 & 14 & 3118\\
4 & 3 & 8 & 53 &       &         &       &       \end{tabular}
\vskip 8pt
Based on empirical evidence, we conjecture that the upper bound in~(5) can be reduced to
\begin{equation}
M<a<\frac{(M+1)(M+2)}3.
\end{equation}
\vskip 24pt\noindent{\bf 4. Large examples of $\boldsymbol{2^kpq}$}
\vskip 12pt\noindent
With $n=2^kpq$ as in \S3, we remark that~$n$ is weird if $a=M+1=2^{k+1}$. For, in this case we have
$$\sum_{j=0}^k2^j<a<p,$$
and so it is impossible to have $a\in\Sigma(\mathcal{A}_n)$. In this case, (2) becomes
\begin{equation}
\frac{p-M}2\cdot\frac{q-M}2=2^k(2^k-1).
\end{equation}
Factoring $2^k-1=uv$ will produce a weird number if~$p$ and~$q$ are both prime, where
\begin{alignat}{2}
p&=M+2^{i+1}u&\;=2^{i+1}(2^{k-i}+u)-1,\\
q&=M+2^{k-i+1}v&\;=2^{k-i+1}(2^{i}+v)-1,\end{alignat}
for some~$i$, $0\le i\le k$.

For speed and simplicity, we factored $2^k-1$ in~(7) algebraically: we considered $k=24j$ and obtained
$$2^k-1=2^{24j}-1=\prod_{d\mid24}\Phi_d(2^j)=uv,$$
where $\Phi_m(x)$ denotes the cyclotomic polynomial of order~$m$ evaluated at~$x$. Then the factors~$u$ and~$v$ were generated by taking all possible subsets
$$A\cup B=\{d:d\mid 24\},\qquad A\cap B=\varnothing,$$
such that
$$u=\prod_{d\in A}\Phi_d(2^j),\qquad v=\prod_{d\in B}\Phi_d(2^j).$$
We also choose all possible $0\le i\le k$ to produce~$p$ and~$q$ in~(8), (9). If such a pair~$p$, $q$ are both prime then $2^kpq$ is weird.

To test $p$ in (8) and $q$ in (9) for primality, we first ran an initial screen comprising pseudoprime tests for the bases~2, 3, 5, and~7. Assuming both~$p$ and~$q$ survive the screen, we then applied Lucas sequences to try and prove primality for both. Recall that $U_n$ and $V_n$ are the Lucas sequences with parameters $P$, $Q$ if for all $n\ge0$
$$U_n=\frac{\alpha^n-\beta^n}{\alpha-\beta},\qquad V_n=\alpha^n+\beta^n,$$
where $\alpha$ and $\beta$ are the real (or complex) numbers determined by $\alpha+\beta=P$, $\alpha\beta=Q$. In all our tests we take $Q=-1$, and $\Delta=P^2+4$. We were able to produce several prime pairs~$p$, $q$ by applying either Theorem~4.2.3 or Theorem~4.2.8, found in Chapter~4, Section~2, Crandall and Pomerance~[2]. We rewrite both theorems here (in their given order) with less generality than they're given in the text~[2], in accordance with our specific needs.
\vskip 12pt\noindent{\bf Proposition 1.}  {\it Suppose $\left(\frac{\Delta}{n}\right)=-1$, $Q=-1$, and~$n$ is an odd positive integer. Suppose $F=2^k$ for some $k>1$, $F\mid n+1$, and that the Lucas sequence $U_n$ has the property
$$U_{n+1}\equiv0\pmod{n},\qquad\gcd(U_{(n+1)/2},n)=1.$$
Then every prime $p$ dividing~$n$ satisfies $p\equiv\left(\frac{\Delta}{n}\right)\pmod{F}$. In particular, if $F>\sqrt{n}+1$, then~$n$ is prime.}
\vskip 12pt\noindent
Note that $\left(\frac{\Delta}{n}\right)$ denotes the Jacobi symbol. 
\vskip 12pt\noindent{\bf Proposition 2.}  {\it Suppose $\left(\frac{\Delta}{n}\right)=-1$, $Q=-1$, and~$n$ is an odd positive integer. Suppose $F=2^k$ for some $k>1$, and $F\mid n+1$ where $F>n^{1/3}+1$. Suppose further that the Lucas sequence $U_n$ has the property
$$U_{n+1}\equiv0\pmod{n},\qquad \gcd(U_{(n+1)/2},n)=1.$$
Write $n+1=FR$, and write~$R$ in base~$F$ so that $R=r_1F+r_0$, $0\le r_i\le F-1$. Then~$n$ is prime if and only if neither $x^2+r_0x-r_1$ nor $x^2+(r_0+F)-r_1-1$ has a positive integral root.}
\vskip 12pt

In~(8), (9), both $p+1$ and $q+1$ have the form $FR$, where $F$ is a power of~2 and~$R$ is odd. We calculate both $U_R$ and $V_R$ by repeated use of the identities
$$U_{2m}=U_mV_m,\qquad V_{2m}=V_m^2-2Q^m,$$
and
$$2U_{m+1}=PU_m+V_m,\qquad 2V_{m+1}=PV_m+\Delta U_m.$$
Since we only use $Q=-1$,  we have $V_{2m}=V_m^2-2$ when~$m$ is even, and $V_{2m}=V_m^2+2$ when~$m$ is odd. 

After computing $U_R$ and $V_R$, we need only compute $V_{FR/2}$. For, we have $\gcd(V_{FR/2},U_{FR/2})=1$ or~2, and $U_{FR}=U_{FR/2}V_{FR/2}$. Since~$p$ and~$q$ are odd, the condition $V_{FR/2}\equiv0\pmod{p}$ (respectively $V_{FR/2}\equiv0\pmod{q}$) implies that $p\mid U_{FR}$, $\gcd(U_{FR/2},p)=1$ (respectively $q\mid U_{FR}$, $\gcd(U_{FR/2},q)=1$). Then $V_{FR/2}$ is computed by applying first $V_{2R}=V_R^2+2$, and then applying\break$V_{2^kR}=V_{2^{k-1}R}^2-2$ for all $k\ge2$. 

We found 66 such primitive weird numbers $2^kpq$, ranging in digit length from 219 to 2077. We list them in Tables~42--48, after the references. For ease in reading these tables, the primes~$p$ and~$q$ are given in the form of~(8), (9), with~$u$ and~$v$ expressed as products of cyclotomic polynomials $\Phi_d(2^j)$. Please bear in mind here that we shall use the symbol $\Phi_d$ to represent $\Phi_d(2^j)$, and~$j$ will be given. 

Recall that
\begin{alignat*}{2}
\Phi_1(x)&=x-1,\qquad&\qquad \Phi_6(x)&=x^2-x+1,\\
\Phi_2(x)&=x+1,\qquad&\qquad \Phi_8(x)&=x^4+1,\\
\Phi_3(x)&=x^2+x+1,\qquad&\qquad \Phi_{12}(x)&=x^4-x^2+1,\\
\Phi_4(x)&=x^2+1,\qquad&\qquad \Phi_{24}(x)&=x^8-x^4+1.\end{alignat*}
With each weird number listed, we indicate the parameters~$P$ and~$\Delta$ of the particular Lucas sequences used to prove the primality of~$p$ and~$q$ respectively. We also write either ``Thm. 4.2.3'' or ``Thm. 4.2.8'' to indicate which of the two theorems by Crandall and Pomerance~[2] were applied in each of the primality proofs for~$p$ and~$q$. We also indicate the digit length of~$p$, $q$, and $n=2^kpq$. 
\vskip 24pt
\noindent{\bf References}
\vskip 12pt
\begin{description}
\item{[1]}\quad S. Benkoski and P. Erd\H{o}s, ``On weird and pseudoperfect numbers,'' {\it Mathematics of Computation\/}, {\bf 28}, No. 126 (april 1974), 617--623. 
\item{[2]}\quad R. Crandall and C. Pomerance, {\it Prime Numbers: A Computational Perspective}, Springer-Verlag (New York), 2001.
\item{[3]}\quad S. Kravitz, ``A search for large weird numbers,'' {\it Journal of Recreational Mathematics\/}, {\bf 9}, No. 2, (1976--77), 82--85.
\item{[4]}\quad S. Pajunen, ``On primitive weird numbers,'' {\it MRFS: A Collection of Manuscripts Related to the Fibonacci Sequence (18$^{th}$ Anniversary Volume)}, V. Hoggat, ed., Fibonacci Association, (1980), 162--166. 
\end{description}
\vskip 24pt\noindent
{\sc Towers Condominiums B-16}

\noindent
{\sc 3600 Contant}

\noindent
{\sc St. Thomas, VI\quad 00802}

\noindent
{\sc USA}

\noindent
{\sc Email:}\quad{\tt diannuc@gmail.com}
\vfill\eject
%
%
%
 
\noindent{\bf Table 1.}\quad All primitive weird numbers of form $2^kpq$, $1\le k\le8$.
\vskip 4pt
\centerline{
\begin{tabular}{|lllll|}\hline
       &                           &                         &                         &                         \\
$k=1$  &$2\cdot5\cdot7$            &                         &  $\qquad k=2$                  & $2^2\cdot11\cdot19$     \\
       &                           &                         &                         &                         \\
$k=3$  &$2^3\cdot19\cdot61$        &$2^3\cdot17\cdot127$     &$2^3\cdot19\cdot71$      &$2^3\cdot23\cdot43$      \\
       &$2^3\cdot29\cdot31$        &                         &                         &                         \\
       &                           &                         &                         &                         \\
$k=4$  &$2^4\cdot43\cdot107$       &$2^4\cdot37\cdot191$     &$2^4\cdot41\cdot127$     &                         \\
       &                           &                         &                         &                         \\
\hline
\end{tabular}}
\vfill\eject
%
%
%
 
\noindent{\bf Table 2.}\quad All primitive weird numbers of form $2^9pq$.
\vskip 6pt
\centerline{\begin{tabular}{|llll|}\hline
                            &                               &                               &                             \\
$2^9\cdot1031\cdot121013$   &   $2^9\cdot1031\cdot121271$   &   $2^9\cdot1109\cdot12211 $   &   $2^9\cdot1031\cdot122173$ \\ 
$2^9\cdot1033\cdot99581 $   &   $2^9\cdot1033\cdot99787 $   &   $2^9\cdot1031\cdot124493$   &   $2^9\cdot1039\cdot62989 $ \\  
$2^9\cdot1039\cdot63443 $   &   $2^9\cdot1033\cdot101027$   &   $2^9\cdot1039\cdot63703 $   &   $2^9\cdot1049\cdot39821 $ \\  
$2^9\cdot1033\cdot101957$   &   $2^9\cdot1039\cdot64223 $   &   $2^9\cdot1049\cdot39983       $   &   $2^9\cdot1049\cdot40063 $ \\  
$2^9\cdot2011\cdot2053              $   &   $2^9\cdot1051\cdot7423 $   &   $2^9\cdot1033\cdot103093$   &   $2^9\cdot1061\cdot27893 $ \\
$2^9\cdot1039\cdot64937 $   &   $2^9\cdot1051\cdot37573 $   &   $2^9\cdot2017\cdot2053           $   &   $2^9\cdot1223\cdot6151              $ \\
$2^9\cdot1063\cdot26681           $   &   $2^9\cdot1049\cdot40507 $   &   $2^9\cdot2011\cdot2063                 $   &   $2^9\cdot1223\cdot6163              $ \\
\hline
\end{tabular}}

\vfill\eject
%
%
%
 
\noindent{\bf Table 3.}\quad All primitive weird numbers of form $2^{10}pq$. (Part 1)
\vskip 6pt
\centerline{\begin{tabular}{|llll|}\hline
 & & & \\
$2^{10}\!\cdot\!2053\!\cdot\!605977$  &         $2^{10}\!\cdot\! 2053\!\cdot\!607003$   &      $2^{10}\!\cdot\!2053\!\cdot\!614531$   &      $2^{10}\!\cdot\!2053\!\cdot\!615557$       \\
$2^{10}\!\cdot\!2053\!\cdot\! 618637$   &       $2^{10}\!\cdot\!2053\!\cdot\!618979 $  &       $2^{10}\!\cdot\!2053\!\cdot\!627191$    &     $2^{10}\!\cdot\!2053\!\cdot\!630613$      \\       
$2^{10}\!\cdot\!2053\!\cdot\!   639167$   &     $2^{10}\!\cdot\!2053\!\cdot\! 639851$   &      $2^{10}\!\cdot\!2053\!\cdot\! 640193$  &      $2^{10}\!\cdot\!2053 \!\cdot\! 642931$    \\   
$2^{10}\!\cdot\!2053 \!\cdot\!  643273$      &  $2^{10}\!\cdot\!2053\!\cdot\!643957$ &         $2^{10}\!\cdot\!2053\!\cdot\!651143$  &       $2^{10}\!\cdot\!2053\!\cdot\!  676463$   \\    
$2^{10}\!\cdot\!2053 \!\cdot\!  677147$     &   $2^{10}\!\cdot\!2053\!\cdot\!680569$  &        $2^{10}\!\cdot\!2053\!\cdot\!681253$ &        $2^{10}\!\cdot\! 2053 \!\cdot\! 688097$\\ 
$2^{10}\!\cdot\!2053\!\cdot\!   693571$      &  $2^{10}\!\cdot\!2063\!\cdot\!249107$  &        $2^{10}\!\cdot\!2063\!\cdot\!250267$   &      $2^{10}\!\cdot\!2063  \!\cdot\!260453$     \\ 
$2^{10}\!\cdot\!2063\!\cdot\!   260969$     &   $2^{10}\!\cdot\!2063\!\cdot\!262387$   &       $2^{10}\!\cdot\!2063\!\cdot\!263677$  &       $2^{10}\!\cdot\!2069 \!\cdot\! 184517 $   \\
\hline\end{tabular}}
\vfill\eject

 %
%
%
 
\noindent{\bf Table 4.}\quad All primitive weird numbers of form $2^{10}pq$. (Part 2)
\vskip 6pt
\centerline{\begin{tabular}{|llll|}\hline
 & & & \\
$2^{10}\!\cdot\!3299 \!\cdot\! 5393$    &            $2^{10}\!\cdot\!3469\!\cdot\!  4993$    &           $2^{10}\!\cdot\!3511\!\cdot\!  4909$    &           $2^{10}\!\cdot\!3607 \!\cdot\! 4733  $ \\           
$2^{10}\!\cdot\!3673\!\cdot\!  4621$    &            $2^{10}\!\cdot\!3739\!\cdot\!  4523$    &           $2^{10}\!\cdot\!3767 \!\cdot\! 4483$    &           $2^{10}\!\cdot\!3797\!\cdot\!  4441  $ \\           
$2^{10}\!\cdot\!3821\!\cdot\!  4409 $    &           $2^{10}\!\cdot\!3881\!\cdot\!  4327$    &           $2^{10}\!\cdot\!3907 \!\cdot\! 4297$    &           $2^{10}\!\cdot\!3917  \!\cdot\!4283$ \\   
$2^{10}\!\cdot\!3919\!\cdot\!  4283$    &            $2^{10}\!\cdot\!3929\!\cdot\!  4271$    &           $2^{10}\!\cdot\!3929 \!\cdot\! 4273 $    &          $2^{10}\!\cdot\!3931 \!\cdot\! 4271  $ \\           
$2^{10}\!\cdot\!3947\!\cdot\!  4241$    &            $2^{10}\!\cdot\!3967\!\cdot\!  4229$    &           $2^{10}\!\cdot\!4013\!\cdot\!  4157$    &           $2^{10}\!\cdot\!4013 \!\cdot\! 4159 $ \\           
$2^{10}\!\cdot\!4019\!\cdot\!  4157$    &            $2^{10}\!\cdot\!4019 \!\cdot\! 4159 $    &          $2^{10}\!\cdot\!4021\!\cdot\!  4157$    &           $2^{10}\!\cdot\!4021\!\cdot\!  4159   $ \\          
$2^{10}\!\cdot\!4027\!\cdot\!  4159$    &            $2^{10}\!\cdot\!4049 \!\cdot\! 4099 $    &          $2^{10}\!\cdot\!4049 \!\cdot\! 4111$    &           $2^{10}\!\cdot\!4049\!\cdot\!  4127  $ \\
\hline\end{tabular}}

\vskip 24pt

\noindent{\bf Table 5.}\quad All primitive weird numbers of form $2^{11}pq$. (Part 1)
\vskip 6pt
\centerline{\begin{tabular}{|llll|}\hline
 & & & \\ 
$2^{11}\!\cdot\!4099 \!\cdot\! 3368353$  &   $2^{11}\!\cdot\!4099\!\cdot\!  3387823  $  &  $2^{11}\!\cdot\!4099\!\cdot\!  3417541$  &    $2^{11}\!\cdot\!4099\!\cdot\!  3448283$\\    
$2^{11}\!\cdot\!4099\!\cdot\!  3486199$  &   $2^{11}\!\cdot\!4099\!\cdot\!  3497471 $  &   $2^{11}\!\cdot\!4099 \!\cdot\! 3507719 $  &   $2^{11}\!\cdot\!4099\!\cdot\!  3515917    $\\
$2^{11}\!\cdot\!4099\!\cdot\!  3532313$  &   $2^{11}\!\cdot\!4099 \!\cdot\! 3589699 $  &   $2^{11}\!\cdot\!4099\!\cdot\!  3595847$  &    $2^{11}\!\cdot\!4099 \!\cdot\! 3650159    $\\
$2^{11}\!\cdot\!4099 \!\cdot\! 3674753$  &   $2^{11}\!\cdot\!4099 \!\cdot\! 3677827$  &  $2^{11}\!\cdot\! 4099 \!\cdot\! 3679877$  &   $2^{11}\!\cdot\!4099 \!\cdot\! 3712669    $\\
$2^{11}\!\cdot\!4099 \!\cdot\! 3732139 $ &  $2^{11}\!\cdot\! 4099 \!\cdot\! 3751609 $  &   $2^{11}\!\cdot\!4099 \!\cdot\! 3756733$  &    $2^{11}\!\cdot\!4099\!\cdot\!  3773129$\\
$2^{11}\!\cdot\!4099 \!\cdot\! 3778253$  &     $2^{11}\!\cdot\!4099 \!\cdot\! 3792599$  &    $2^{11}\!\cdot\!4099 \!\cdot\! 3797723 $  &   $2^{11}\!\cdot\!4099\!\cdot\!  3817193  $\\  
$2^{11}\!\cdot\!4099 \!\cdot\! 3822317$  &     $2^{11}\!\cdot\!4099 \!\cdot\! 3827441$  &    $2^{11}\!\cdot\!4099 \!\cdot\! 3841787 $  &  $2^{11}\!\cdot\! 4099 \!\cdot\! 3874579 $\\
\hline\end{tabular}}

\noindent{\bf Table 6.}\quad All primitive weird numbers of form $2^{11}pq$. (Part 2)
\vskip 6pt
\centerline{\begin{tabular}{|llll|}\hline
 & & & \\ 
$2^{11}\!\cdot\!4157\!\cdot\!  272417 $  &      $2^{11}\!\cdot\! 4157 \!\cdot\! 273289   $  &    $2^{11}\!\cdot\!4157 \!\cdot\! 274093  $  &    $2^{11}\!\cdot\! 4159\!\cdot\!  264031$\\        
$2^{11}\!\cdot\!4159 \!\cdot\! 264811$  &        $2^{11}\!\cdot\!4159\!\cdot\!  266111   $  &    $2^{11}\!\cdot\!4177 \!\cdot\! 206303  $  &    $2^{11}\!\cdot\! 4177 \!\cdot\! 207169   $\\     
$2^{11}\!\cdot\!4177 \!\cdot\! 208391 $  &       $2^{11}\!\cdot\!4177 \!\cdot\! 208493   $  &    $2^{11}\!\cdot\!4201 \!\cdot\! 161461 $  &     $2^{11}\!\cdot\! 4201 \!\cdot\! 161659    $\\ 
$2^{11}\!\cdot\!4201 \!\cdot\! 162293 $  &      $2^{11}\!\cdot\!4211 \!\cdot\! 147457   $  &    $2^{11}\!\cdot\!4211  \!\cdot\!148147  $  &     $2^{11}\!\cdot\!4217 \!\cdot\! 140681 $\\       
$2^{11}\!\cdot\!4217\!\cdot\!  141131$  &        $2^{11}\!\cdot\!4217 \!\cdot\! 141269 $  &      $2^{11}\!\cdot\!4217 \!\cdot\! 141511 $  &     $2^{11}\!\cdot\! 4229  \!\cdot\!128321$\\  
$2^{11}\!\cdot\!4231\!\cdot\!  126743 $  &       $2^{11}\!\cdot\!4231\!\cdot\!  126961  $  &     $2^{11}\!\cdot\!4231 \!\cdot\! 127241  $  &    $2^{11}\!\cdot\! 4241\!\cdot\!  118399     $\\   
$2^{11}\!\cdot\!4241\!\cdot\!  118457 $  &       $2^{11}\!\cdot\!4243\!\cdot\!  116797   $  &    $2^{11}\!\cdot\!4243 \!\cdot\! 116911  $  &     $2^{11}\!\cdot\!4243\!\cdot\!  116969  $\\
\hline\end{tabular}}

\noindent{\bf Table 7.}\quad All primitive weird numbers of form $2^{11}pq$. (Part 3)
\vskip 6pt
\centerline{\begin{tabular}{|llll|}\hline
 & & & \\ 
$2^{11}\!\cdot\!5783\!\cdot\!   14029 $  &         $2^{11}\!\cdot\!5801\!\cdot\!   13921  $  &      $2^{11}\!\cdot\! 5821\!\cdot\!   13807  $  &    $2^{11}\!\cdot\! 5897\!\cdot\!   13397  $\\       
$2^{11}\!\cdot\!5927 \!\cdot\!  13241 $  &         $2^{11}\!\cdot\!5981 \!\cdot\!  12979   $  &     $2^{11}\!\cdot\! 5981\!\cdot\!   12983   $  &     $2^{11}\!\cdot\! 6043 \!\cdot\!  12703  $\\       
$2^{11}\!\cdot\!6053\!\cdot\!  12659 $  &        $2^{11}\!\cdot\! 6067\!\cdot\!   12589  $  &     $2^{11}\!\cdot\!  6073 \!\cdot\!  12569   $  &    $2^{11}\!\cdot\!  6079\!\cdot\!   12547$\\ 
$2^{11}\!\cdot\!6091 \!\cdot\! 12487  $  &        $2^{11}\!\cdot\!6101 \!\cdot\! 12451  $  &      $2^{11}\!\cdot\! 6113\!\cdot\!  12289  $  &     $2^{11}\!\cdot\!  6113 \!\cdot\! 12301  $\\       
$2^{11}\!\cdot\!6113 \!\cdot\!12377 $  &         $2^{11}\!\cdot\!6113\!\cdot\!  12401 $  &       $2^{11}\!\cdot\! 6121\!\cdot\!  12323  $  &     $2^{11}\!\cdot\!  6121  \!\cdot\!12329     $\\    
$2^{11}\!\cdot\!6121\!\cdot\!  12347$  &          $2^{11}\!\cdot\!6131 \!\cdot\! 12289 $  &      $2^{11}\!\cdot\!  6131\!\cdot\!  12301  $  &     $2^{11}\!\cdot\!  6133 \!\cdot\! 12323    $\\     
$2^{11}\!\cdot\!6143  \!\cdot\!12211  $  &       $2^{11}\!\cdot\!6143 \!\cdot\! 12241  $  &      $2^{11}\!\cdot\! 6143 \!\cdot\! 12253   $  &     $2^{11}\!\cdot\! 6143 \!\cdot\! 12277   $\\
\hline\end{tabular}}

\vfill\eject
\noindent{\bf Table 8.}\quad All primitive weird numbers of form $2^{12}pq$. (Part 1)
\vskip 6pt
\centerline{\begin{tabular}{|llll|}\hline
 & & & \\ 
$2^{12}\!\cdot\!8209\!\cdot\!   3541271$  &    $2^{12}\!\cdot\!8209\!\cdot\!    3545831$  &  $2^{12}\!\cdot\!8209\!\cdot\!   3556321$  &  $2^{12}\!\cdot\!8209\!\cdot\!   3557689$\\      
$2^{12}\!\cdot\!8209\!\cdot\!   3606487$  &    $2^{12}\!\cdot\!8209\!\cdot\!    3608767$  &  $2^{12}\!\cdot\!8209\!\cdot\!   3610591$  &  $2^{12}\!\cdot\!8209\!\cdot\!   3637043$\\      
$2^{12}\!\cdot\!8209\!\cdot\!   3639323$  &    $2^{12}\!\cdot\!8209\!\cdot\!    3641147$  &  $2^{12}\!\cdot\!8209\!\cdot\!   3642971$  &  $2^{12}\!\cdot\!8209\!\cdot\!   3643883$\\      
$2^{12}\!\cdot\!8209\!\cdot\!   3644339$  &    $2^{12}\!\cdot\!8209\!\cdot\!    3654829$  &  $2^{12}\!\cdot\!8209\!\cdot\!   3655741$  &  $2^{12}\!\cdot\!8209\!\cdot\!   3658021$\\      
$2^{12}\!\cdot\!8209\!\cdot\!   3658477$  &    $2^{12}\!\cdot\!8209\!\cdot\!    3687209$  &  $2^{12}\!\cdot\!8209\!\cdot\!   3690857$  &  $2^{12}\!\cdot\!8209\!\cdot\!   3691313$\\ 
$2^{12}\!\cdot\!8209\!\cdot\!   3692681$  &    $2^{12}\!\cdot\!8209\!\cdot\!    3703171$  &  $2^{12}\!\cdot\!8209\!\cdot\!   3705451$  &  $2^{12}\!\cdot\!8209\!\cdot\!   3706819$\\      
$2^{12}\!\cdot\!8209\!\cdot\!   3708643$  &    $2^{12}\!\cdot\!8219\!\cdot\!    2329499$  &  $2^{12}\!\cdot\!8219\!\cdot\!   2339773$  &  $2^{12}\!\cdot\!8219\!\cdot\!   2343589$\\
\hline\end{tabular}}

\vfill\eject
\noindent{\bf Table 9.}\quad All primitive weird numbers of form $2^{12}pq$. (Part 2)
\vskip 6pt
\centerline{\begin{tabular}{|llll|}\hline
 & & & \\ 
 $2^{12}\!\cdot\!8377 \!\cdot\!   368047 $  &      $2^{12}\!\cdot\!8377\!\cdot\!    368227$  &      $2^{12}\!\cdot\!8377  \!\cdot\!  368633 $  &     $2^{12}\!\cdot\!8387 \!\cdot\!   349387 $\\     
$2^{12}\!\cdot\!8387\!\cdot\!    349729 $  &      $2^{12}\!\cdot\!8389 \!\cdot\!   345431$  &      $2^{12}\!\cdot\!8389  \!\cdot\!  345601 $  &     $2^{12}\!\cdot\!8389  \!\cdot\!  345643 $\\  
 $2^{12}\!\cdot\!8389 \!\cdot\!   346321 $  &      $2^{12}\!\cdot\!8389 \!\cdot\!   347041  $  &   $2^{12}\!\cdot\! 8419\!\cdot\!    301237$  &    $2^{12}\!\cdot\! 8419  \!\cdot\!  301459$\\  
$2^{12}\!\cdot\!8419\!\cdot\!    301643 $  &      $2^{12}\!\cdot\!8419 \!\cdot\!   302123 $  &    $2^{12}\!\cdot\! 8423 \!\cdot\!   296729 $  &     $2^{12}\!\cdot\!8423 \!\cdot\!   296801  $\\    
$2^{12}\!\cdot\!8423 \!\cdot\!   296983 $  &      $2^{12}\!\cdot\!8423 \!\cdot\!   297019 $  &     $2^{12}\!\cdot\!8429 \!\cdot\!   289171$  &     $2^{12}\!\cdot\! 8429  \!\cdot\!  290021 $\\     
$2^{12}\!\cdot\!8431 \!\cdot\!   287251 $  &      $2^{12}\!\cdot\!8431\!\cdot\!    287321 $  &    $2^{12}\!\cdot\! 8443 \!\cdot\!   270209 $  &    $2^{12}\!\cdot\! 8443  \!\cdot\!  273727  $\\    
$2^{12}\!\cdot\!8443 \!\cdot\!   273827 $  &      $2^{12}\!\cdot\!8447\!\cdot\!    269281 $  &    $2^{12}\!\cdot\!8447  \!\cdot\!  269413 $  &     $2^{12}\!\cdot\!8447  \!\cdot\!  270073   $\\
\hline\end{tabular}}

\vfill\eject
\noindent{\bf Table 10.}\quad All primitive weird numbers of form $2^{12}pq$. (Part 3)
\vskip 6pt
\centerline{\begin{tabular}{|llll|}\hline
 & & & \\ 
$2^{12}\!\cdot\!9151 \!\cdot\!  78079  $  &      $2^{12}\!\cdot\! 9161 \!\cdot\!  77339  $  &     $2^{12}\!\cdot\! 9199 \!\cdot\!  74687 $  &       $2^{12}\!\cdot\!9209  \!\cdot\! 73681$\\ 
$2^{12}\!\cdot\!9209  \!\cdot\! 73699 $  &       $2^{12}\!\cdot\!9209  \!\cdot\! 73771 $  &      $2^{12}\!\cdot\!9209 \!\cdot\!  73961  $  &      $2^{12}\!\cdot\!9209  \!\cdot\! 74051   $\\     
$2^{12}\!\cdot\!9227 \!\cdot\!  72907  $  &       $2^{12}\!\cdot\!9241 \!\cdot\!  72053  $  &      $2^{12}\!\cdot\!9257 \!\cdot\!  71129  $  &     $2^{12}\!\cdot\! 9281  \!\cdot\! 69709  $\\      
$2^{12}\!\cdot\!9293 \!\cdot\!  69031 $  &       $2^{12}\!\cdot\! 9293 \!\cdot\!  69073 $  &       $2^{12}\!\cdot\!9311\!\cdot\!   68053 $  &      $2^{12}\!\cdot\! 9323  \!\cdot\! 67427    $\\    
$2^{12}\!\cdot\!9341 \!\cdot\!  66491  $  &       $2^{12}\!\cdot\!9341 \!\cdot\!  66499  $  &      $2^{12}\!\cdot\!9343 \!\cdot\!  66431 $  &      $2^{12}\!\cdot\! 9349  \!\cdot\! 65539   $\\     
$2^{12}\!\cdot\!9371\!\cdot\!   65033  $  &       $2^{12}\!\cdot\!9397 \!\cdot\!  63799  $  &     $2^{12}\!\cdot\! 9397 \!\cdot\!  63823  $  &     $2^{12}\!\cdot\! 9413  \!\cdot\! 63079$\\ 
$2^{12}\!\cdot\!9431\!\cdot\!   62297$  &         $2^{12}\!\cdot\!9463 \!\cdot\!  60899$  &        $2^{12}\!\cdot\!9467\!\cdot\!   60763 $  &      $2^{12}\!\cdot\! 9479  \!\cdot\! 60251 $\\
\hline\end{tabular}}

\vfill\eject
\noindent{\bf Table 11.}\quad All primitive weird numbers of form $2^{12}pq$. (Part 4)
\vskip 6pt
\centerline{\begin{tabular}{|llll|}\hline
 & & & \\ 
 $2^{12}\!\cdot\!11161\!\cdot\!  30781$  &        $2^{12}\!\cdot\!11279  \!\cdot\!29917 $  &      $2^{12}\!\cdot\!11383 \!\cdot\! 29209 $  &      $2^{12}\!\cdot\!11399 \!\cdot\! 29101    $\\   
$2^{12}\!\cdot\!11411\!\cdot\!  29023$  &        $2^{12}\!\cdot\!11411\!\cdot\!  29027  $  &    $2^{12}\!\cdot\!11423  \!\cdot\!28949  $  &     $2^{12}\!\cdot\!11447 \!\cdot\! 28793$\\
$2^{12}\!\cdot\!11483\!\cdot\!   28571   $  &     $2^{12}\!\cdot\!11519 \!\cdot\!  28351  $  &     $2^{12}\!\cdot\!11587  \!\cdot\! 27943   $  &    $2^{12}\!\cdot\!11587\!\cdot\!   27947       $\\ 
 $2^{12}\!\cdot\!11621 \!\cdot\!  27751 $  &       $2^{12}\!\cdot\!11677 \!\cdot\!  27437$  &       $2^{12}\!\cdot\!11689 \!\cdot\!  27367   $  &    $2^{12}\!\cdot\!11731 \!\cdot\!  27143       $\\    
$2^{12}\!\cdot\!11777 \!\cdot\!  26893  $  &      $2^{12}\!\cdot\!11813 \!\cdot\!  26711   $  &    $2^{12}\!\cdot\!11863  \!\cdot\! 26459 $  &     $2^{12}\!\cdot\! 11971\!\cdot\!   25933       $\\   
$2^{12}\!\cdot\!12041\!\cdot\!   25601  $  &      $2^{12}\!\cdot\!12097\!\cdot\!   25367  $  &     $2^{12}\!\cdot\!12161\!\cdot\!   25087 $  &      $2^{12}\!\cdot\!12163  \!\cdot\! 25057      $\\  
$2^{12}\!\cdot\!12197\!\cdot\!   24917$  &       $2^{12}\!\cdot\!12197 \!\cdot\!  24923  $  &    $2^{12}\!\cdot\!12203 \!\cdot\!  24889 $  &      $2^{12}\!\cdot\!12203  \!\cdot\! 24907    $\\
\hline\end{tabular}}

\vfill\eject
\noindent{\bf Table 12.}\quad All primitive weird numbers of form $2^{12}pq$. (Part 5)
\vskip 6pt
\centerline{\begin{tabular}{|llll|}\hline
 & & & \\ 
$2^{12}\!\cdot\!16193 \!\cdot\!  16567 $  &       $2^{12}\!\cdot\!16193 \!\cdot\!  16573  $  &     $2^{12}\!\cdot\!16223 \!\cdot\!  16519   $  &    $2^{12}\!\cdot\!16229 \!\cdot\!  16519       $\\     
$2^{12}\!\cdot\!16229 \!\cdot\!  16529 $  &       $2^{12}\!\cdot\!16231  \!\cdot\! 16519  $  &     $2^{12}\!\cdot\!16231 \!\cdot\!  16529  $  &     $2^{12}\!\cdot\!16273  \!\cdot\! 16411    $\\  
$2^{12}\!\cdot\!16273 \!\cdot\!  16417 $  &       $2^{12}\!\cdot\!16273\!\cdot\!   16421$  &       $2^{12}\!\cdot\!16301 \!\cdot\!  16411  $  &     $2^{12}\!\cdot\!16301 \!\cdot\!  16417       $\\     
$2^{12}\!\cdot\!16301 \!\cdot\!  16421$  &        $2^{12}\!\cdot\!16301 \!\cdot\!  16427  $  &     $2^{12}\!\cdot\!16301 \!\cdot\!  16433   $  &    $2^{12}\!\cdot\!16301 \!\cdot\!  16447       $\\     
$2^{12}\!\cdot\!16301\!\cdot\!   16451 $  &       $2^{12}\!\cdot\!16301 \!\cdot\!  16453  $  &     $2^{12}\!\cdot\!16319 \!\cdot\!  16411 $  &      $2^{12}\!\cdot\!16319 \!\cdot\!  16417      $\\      
$2^{12}\!\cdot\!16319 \!\cdot\!  16421$  &        $2^{12}\!\cdot\!16319 \!\cdot\!  16427  $  &     $2^{12}\!\cdot\!16319 \!\cdot\!  16433 $  &      $2^{12}\!\cdot\!16333\!\cdot\!   16411     $\\       
$2^{12}\!\cdot\!16333\!\cdot\!   16417$  &        $2^{12}\!\cdot\!16333 \!\cdot\!  16421$  &       $2^{12}\!\cdot\!16333 \!\cdot\! 16427  $  &     $2^{12}\!\cdot\!16339 \!\cdot\!  16411    $\\  
$2^{12}\!\cdot\!16339\!\cdot\!  16417$  &     $2^{12}\!\cdot\!16339\!\cdot\!   16421 $  &      $2^{12}\!\cdot\!16349 \!\cdot\!  16411$    &      \\  
       &  &   & \\\hline\end{tabular}}

\vskip 12pt

\noindent{\bf Table 13.}\quad All primitive weird numbers of form $2^{13}pq$. (Part 1)
\vskip 6pt
\centerline{\begin{tabular}{|llll|}\hline
 & & & \\  
$2^{13}\!\cdot\!16411\!\cdot\!  9272801 $  &  $2^{13}\!\cdot\!16411  \!\cdot\!9273973$  &   $2^{13}\!\cdot\!16411 \!\cdot\! 9275731 $  &  $2^{13}\!\cdot\!16411 \!\cdot\! 9305623   $\\                              
$2^{13}\!\cdot\!16411 \!\cdot\! 9316759 $  &  $2^{13}\!\cdot\!16411 \!\cdot\! 9319103 $  &  $2^{13}\!\cdot\!16411 \!\cdot\! 9328481 $  &  $2^{13}\!\cdot\!16411 \!\cdot\! 9330239   $\\                                
$2^{13}\!\cdot\!16411\!\cdot\!  9350753 $  & $2^{13}\!\cdot\!16411 \!\cdot\! 9363061  $  & $2^{13}\!\cdot\!16411 \!\cdot\! 9363647$  &   $2^{13}\!\cdot\!16411 \!\cdot\! 9372439    $\\               
$2^{13}\!\cdot\!16411\!\cdot\!  9381817$  & $2^{13}\!\cdot\!16411 \!\cdot\! 9385919$  &   $2^{13}\!\cdot\!16411 \!\cdot\! 9392953$  &   $2^{13}\!\cdot\!16411\!\cdot\!  9416983     $\\              
$2^{13}\!\cdot\!16411\!\cdot\!  9418741 $  & $2^{13}\!\cdot\!16411\!\cdot\!  9419327 $  & $2^{13}\!\cdot\!16411\!\cdot\!  9428119 $  &  $2^{13}\!\cdot\!16411\!\cdot\!  9438083 $\\   
$2^{13}\!\cdot\!16411 \!\cdot\! 9441599 $  & $2^{13}\!\cdot\!16411 \!\cdot\! 9450977 $  &  $2^{13}\!\cdot\!16411\!\cdot\!  9471491$  &   $2^{13}\!\cdot\!16411 \!\cdot\! 9480283  $\\                
$2^{13}\!\cdot\!16411 \!\cdot\! 9503141$  &  $2^{13}\!\cdot\!16411 \!\cdot\! 9504899$  &  $2^{13}\!\cdot\! 16411 \!\cdot\! 9506657$  &  $2^{13}\!\cdot\! 16411 \!\cdot\! 9525413   $\\
\hline\end{tabular}}

\vfill\eject
\noindent{\bf Table 14.}\quad All primitive weird numbers of form $2^{13}pq$. (Part 2)
\vskip 6pt
\centerline{\begin{tabular}{|llll|}\hline
 & & & \\  
$2^{13}\!\cdot\!16433\!\cdot\!  5314103$  &  $2^{13}\!\cdot\!  16433\!\cdot\!  5317061 $  &  $2^{13}\!\cdot\!  16433\!\cdot\!  5326921  $  &  $2^{13}\!\cdot\! 16433 \!\cdot\! 5330207$\\ 
$2^{13}\!\cdot\!16433\!\cdot\!  5340067$  &  $2^{13}\!\cdot\!  16433\!\cdot\!  5342039   $  &  $2^{13}\!\cdot\!16433 \!\cdot\! 5345983 $  &  $2^{13}\!\cdot\!  16433\!\cdot\!  5362087   $\\             
$2^{13}\!\cdot\!16433\!\cdot\!  5363731$  &  $2^{13}\!\cdot\!  16433 \!\cdot\! 5367017 $  &  $2^{13}\!\cdot\!  16433 \!\cdot\! 5368003 $  &  $2^{13}\!\cdot\!  16433\!\cdot\!  5368661  $\\              
$2^{13}\!\cdot\!16433\!\cdot\!  5370961$  &  $2^{13}\!\cdot\!  16433 \!\cdot\! 5372933  $  &  $2^{13}\!\cdot\! 16433 \!\cdot\! 5375891 $  &  $2^{13}\!\cdot\!  16447\!\cdot\!  4150297   $\\             
$2^{13}\!\cdot\!16447 \!\cdot\! 4152353 $  &  $2^{13}\!\cdot\! 16447 \!\cdot\! 4154923  $  &  $2^{13}\!\cdot\! 16447 \!\cdot\! 4158521 $  &  $2^{13}\!\cdot\!  16447 \!\cdot\! 4160063 $\\               
$2^{13}\!\cdot\!16447\!\cdot\!  4162889$  &  $2^{13}\!\cdot\!  16447 \!\cdot\! 4169057  $  &  $2^{13}\!\cdot\! 16447\!\cdot\!  4170599 $  &  $2^{13}\!\cdot\!  16447\!\cdot\!  4172141$\\ 
$2^{13}\!\cdot\!16447 \!\cdot\! 4173683$  &  $2^{13}\!\cdot\!  16447\!\cdot\!  4174711  $  &  $2^{13}\!\cdot\! 16447\!\cdot\!  4176253  $  &  $2^{13}\!\cdot\! 16447\!\cdot\!  4178051 $\\
\hline\end{tabular}}

\vfill\eject
\noindent{\bf Table 15.}\quad All primitive weird numbers of form $2^{13}pq$. (Part 3)
\vskip 6pt
\centerline{\begin{tabular}{|llll|}\hline
 & & & \\  
$2^{13}\!\cdot\!16573\!\cdot\!  1424231$  &  $2^{13}\!\cdot\!  16573  \!\cdot\!1425539  $  &  $2^{13}\!\cdot\! 16573  \!\cdot\!1426063 $  &  $2^{13}\!\cdot\!  16573 \!\cdot\! 1426237 $\\   
$2^{13}\!\cdot\!16573 \!\cdot\! 1426499 $  &  $2^{13}\!\cdot\! 16573  \!\cdot\!1426673  $  &  $2^{13}\!\cdot\! 16573  \!\cdot\!1428767 $  &  $2^{13}\!\cdot\!  16603  \!\cdot\!1232999 $\\               
$2^{13}\!\cdot\!16603 \!\cdot\! 1233301$  &  $2^{13}\!\cdot\!  16603  \!\cdot\!1233527  $  &  $2^{13}\!\cdot\! 16603 \!\cdot\! 1235263  $  &  $2^{13}\!\cdot\! 16603  \!\cdot\!1235791   $\\             
$2^{13}\!\cdot\!16603\!\cdot\!  1236017$  &  $2^{13}\!\cdot\!  16619 \!\cdot\! 1146781 $  &  $2^{13}\!\cdot\!  16619 \!\cdot\! 1149457  $  &  $2^{13}\!\cdot\! 16619 \!\cdot\! 1149527$\\
$2^{13}\!\cdot\!16619 \!\cdot\! 1151147 $  &  $2^{13}\!\cdot\! 16619 \!\cdot\! 1151569  $  &  $2^{13}\!\cdot\! 16631 \!\cdot\!  1094963  $  &  $2^{13}\!\cdot\! 16631  \!\cdot\! 1095433 $\\               
$2^{13}\!\cdot\!16631 \!\cdot\! 1097377 $  &  $2^{13}\!\cdot\! 16633\!\cdot\!  1089461  $  &  $2^{13}\!\cdot\! 16633 \!\cdot\!  1089661  $  &  $2^{13}\!\cdot\! 16649  \!\cdot\! 1021973   $\\             
$2^{13}\!\cdot\!16649 \!\cdot\! 1022849 $  &  $2^{13}\!\cdot\! 16649 \!\cdot\! 1023037 $  &  $2^{13}\!\cdot\! 16651  \!\cdot\! 1014779   $  &  $2^{13}\!\cdot\!16651  \!\cdot\! 1015897$\\
\hline\end{tabular}}

\vfill\eject
\noindent{\bf Table 16.}\quad All primitive weird numbers of form $2^{13}pq$. (Part 4)
\vskip 6pt
\centerline{\begin{tabular}{|llll|}\hline
 & & & \\  
$2^{13}\!\cdot\!17033  \!\cdot\!425903 $  &  $2^{13}\!\cdot\!    17033 \!\cdot\! 428759   $  &  $2^{13}\!\cdot\!  17041 \!\cdot\!  423823   $  &  $2^{13}\!\cdot\! 17047  \!\cdot\! 420193   $\\  
$2^{13}\!\cdot\!17077 \!\cdot\! 402613  $  &  $2^{13}\!\cdot\!   17077\!\cdot\!  402761 $  &  $2^{13}\!\cdot\!    17077  \!\cdot\! 402859  $  &  $2^{13}\!\cdot\!  17093  \!\cdot\! 393187    $\\ 
$2^{13}\!\cdot\!17093 \!\cdot\! 393331 $  &  $2^{13}\!\cdot\!    17093 \!\cdot\! 393403  $  &  $2^{13}\!\cdot\!   17093 \!\cdot\!  393451  $  &  $2^{13}\!\cdot\!  17099\!\cdot\!  391103  $\\   
$2^{13}\!\cdot\!17099\!\cdot\!  391151$  &  $2^{13}\!\cdot\!     17099\!\cdot\!  391199   $  &  $2^{13}\!\cdot\!  17099\!\cdot\!   391247  $  &  $2^{13}\!\cdot\!  17107\!\cdot\!   386963   $\\  
$2^{13}\!\cdot\!17117  \!\cdot\!381937  $  &  $2^{13}\!\cdot\!   17123 \!\cdot\! 376729   $  &  $2^{13}\!\cdot\!  17137 \!\cdot\!  372013 $  &  $2^{13}\!\cdot\!   17137\!\cdot\!   372059 $\\    
$2^{13}\!\cdot\!17137  \!\cdot\!372241 $  &  $2^{13}\!\cdot\!    17137  \!\cdot\!372263   $  &  $2^{13}\!\cdot\! 17159 \!\cdot\!  362107  $  &  $2^{13}\!\cdot\!  17167 \!\cdot\!  358601  $\\    
$2^{13}\!\cdot\!17167\!\cdot\!  358667 $  &  $2^{13}\!\cdot\!    17167  \!\cdot\!358711 $  &  $2^{13}\!\cdot\!    17167  \!\cdot\! 358733   $  &  $2^{13}\!\cdot\! 17189 \!\cdot\!  349303   $\\
\hline\end{tabular}}

\vfill\eject
\noindent{\bf Table 17.}\quad All primitive weird numbers of form $2^{13}pq$. (Part 5)
\vskip 6pt
\centerline{\begin{tabular}{|llll|}\hline
 & & & \\  
$2^{13}\!\cdot\!18199  \!\cdot\!164051 $  &  $2^{13}\!\cdot\!    18199 \!\cdot\! 164071 $  &  $2^{13}\!\cdot\!    18211\!\cdot\!   163151  $  &  $2^{13}\!\cdot\!  18211 \!\cdot\! 163171 $\\ 
$2^{13}\!\cdot\!18211  \!\cdot\!163181$  &  $2^{13}\!\cdot\!     18211\!\cdot\!  163211 $  &  $2^{13}\!\cdot\!    18217\!\cdot\!   162671  $  &  $2^{13}\!\cdot\!  18217 \!\cdot\! 162691 $\\   
$2^{13}\!\cdot\!18217\!\cdot\!  162731$  &  $2^{13}\!\cdot\!    18229\!\cdot\!  161779  $  &  $2^{13}\!\cdot\!   18257 \!\cdot\!  159539   $  &  $2^{13}\!\cdot\! 18269  \!\cdot\!158647  $\\  
$2^{13}\!\cdot\!18269 \!\cdot\! 158657 $  &  $2^{13}\!\cdot\!    18287 \!\cdot\! 157321  $  &  $2^{13}\!\cdot\!   18311\!\cdot\!   155539  $  &  $2^{13}\!\cdot\!  18313 \!\cdot\! 155413$\\    
$2^{13}\!\cdot\!18313\!\cdot\!  155423$  &  $2^{13}\!\cdot\!     18371 \!\cdot\! 151357  $  &  $2^{13}\!\cdot\!   18401 \!\cdot\!  149323 $  &  $2^{13}\!\cdot\!   18401 \!\cdot\! 149341 $\\  
$2^{13}\!\cdot\!18413 \!\cdot\! 148537 $  &  $2^{13}\!\cdot\!    18413 \!\cdot\! 148573   $  &  $2^{13}\!\cdot\!  18427 \!\cdot\!  147551 $  &  $2^{13}\!\cdot\!   18433 \!\cdot\! 147221    $\\  
$2^{13}\!\cdot\!18433  \!\cdot\!147293  $  &  $2^{13}\!\cdot\!   18433 \!\cdot\! 147311  $  &  $2^{13}\!\cdot\!   18439 \!\cdot\!  146857  $  &  $2^{13}\!\cdot\!  18439 \!\cdot\! 146893    $\\
\hline\end{tabular}}

\vfill\eject

\vfill\eject
\noindent{\bf Table 18.}\quad All primitive weird numbers of form $2^{13}pq$. (Part 6)
\vskip 6pt
\centerline{\begin{tabular}{|llll|}\hline
 & & & \\  
$2^{13}\!\cdot\!20483 \!\cdot\! 81737$  &  $2^{13}\!\cdot\!        20483  \!\cdot\!81817  $  &  $2^{13}\!\cdot\!      20483 \!\cdot\!  81847   $  &  $2^{13}\!\cdot\!     20507 \!\cdot\! 81401   $\\    
$2^{13}\!\cdot\!20507\!\cdot\!  81421$  &  $2^{13}\!\cdot\!        20509  \!\cdot\!81409  $  &  $2^{13}\!\cdot\!      20533 \!\cdot\!  80989   $  &  $2^{13}\!\cdot\!     20533 \!\cdot\! 81043     $\\  
$2^{13}\!\cdot\!20543\!\cdot\!  80897$  &  $2^{13}\!\cdot\!        20549  \!\cdot\!80809 $  &  $2^{13}\!\cdot\!       20551\!\cdot\!   80749   $  &  $2^{13}\!\cdot\!     20551\!\cdot\!  80779      $\\ 
$2^{13}\!\cdot\!20593  \!\cdot\!80107  $  &  $2^{13}\!\cdot\!      20693 \!\cdot\! 78643   $  &  $2^{13}\!\cdot\!     20717 \!\cdot\!  78307  $  &  $2^{13}\!\cdot\!   20719\!\cdot\!  78283        $\\ 
$2^{13}\!\cdot\!20743 \!\cdot\! 77933$  &  $2^{13}\!\cdot\!        20749 \!\cdot\! 77849   $  &  $2^{13}\!\cdot\!     20753 \!\cdot\!  77797    $  &  $2^{13}\!\cdot\!    20759  \!\cdot\!77713     $\\  
$2^{13}\!\cdot\!20773\!\cdot\!  77513 $  &  $2^{13}\!\cdot\!       20789 \!\cdot\! 77291  $  &  $2^{13}\!\cdot\!      20807 \!\cdot\!  77047    $  &  $2^{13}\!\cdot\!    20857 \!\cdot\!76369     $\\  
$2^{13}\!\cdot\!20879 \!\cdot\! 76081 $  &  $2^{13}\!\cdot\!       20897 \!\cdot\! 75833   $  &  $2^{13}\!\cdot\!     20959 \!\cdot\!  75037   $  &  $2^{13}\!\cdot\!     20981  \!\cdot\!74747    $\\
\hline\end{tabular}}

\vfill\eject

\vfill\eject
\noindent{\bf Table 19.}\quad All primitive weird numbers of form $2^{13}pq$. (Part 7)
\vskip 6pt
\centerline{\begin{tabular}{|llll|}\hline& & & \\  
$2^{13}\!\cdot\!23537\!\cdot\!  53897 $  &  $2^{13}\!\cdot\!      23539  \!\cdot\!53887  $  &  $2^{13}\!\cdot\!      23593 \!\cdot\!  53609   $  &  $2^{13}\!\cdot\!     23609\!\cdot\!  53527  $\\  
$2^{13}\!\cdot\!23689 \!\cdot\! 53113   $  &  $2^{13}\!\cdot\!     23689  \!\cdot\!53117  $  &  $2^{13}\!\cdot\!      23857  \!\cdot\! 52291   $  &  $2^{13}\!\cdot\!     23879 \!\cdot\! 52189    $\\   
$2^{13}\!\cdot\!23887\!\cdot\!  52147   $  &  $2^{13}\!\cdot\!     23957  \!\cdot\!51817 $  &  $2^{13}\!\cdot\!       24001 \!\cdot\!  51599   $  &  $2^{13}\!\cdot\!     24007 \!\cdot\! 51581    $\\   
$2^{13}\!\cdot\!24023 \!\cdot\! 51511   $  &  $2^{13}\!\cdot\!     24043 \!\cdot\! 51419  $  &  $2^{13}\!\cdot\!      24121 \!\cdot\!  51059   $  &  $2^{13}\!\cdot\!     24169 \!\cdot\! 50839   $\\    
$2^{13}\!\cdot\!24203 \!\cdot\! 50683   $  &  $2^{13}\!\cdot\!     24229 \!\cdot\! 50591  $  &  $2^{13}\!\cdot\!      24247  \!\cdot\! 50513  $  &  $2^{13}\!\cdot\!      24281 \!\cdot\! 50363     $\\  
$2^{13}\!\cdot\!24379 \!\cdot\! 49937   $  &  $2^{13}\!\cdot\!     24407\!\cdot\!  49811 $  &  $2^{13}\!\cdot\!       24407\!\cdot\!   49823  $  &  $2^{13}\!\cdot\!      24413  \!\cdot\!49783     $\\   
$2^{13}\!\cdot\!24413 \!\cdot\! 49789 $  &  $2^{13}\!\cdot\!       24413 \!\cdot\! 49801  $  &  $2^{13}\!\cdot\!      24413 \!\cdot\!  49807    $  &  $2^{13}\!\cdot\!    24419 \!\cdot\! 49783   $\\
\hline\end{tabular}}

\vfill\eject

\vfill\eject
\noindent{\bf Table 20.}\quad All primitive weird numbers of form $2^{13}pq$. (Part 8)
\vskip 6pt
\centerline{\begin{tabular}{|llll|}\hline& & & \\  
$2^{13}\!\cdot\!29437 \!\cdot\! 36943     $  &  $2^{13}\!\cdot\!   29473 \!\cdot\! 36887   $  &  $2^{13}\!\cdot\!     29483 \!\cdot\!  36871  $  &  $2^{13}\!\cdot\!      29537  \!\cdot\!36787   $\\ 
$2^{13}\!\cdot\!29599 \!\cdot\! 36691   $  &  $2^{13}\!\cdot\!     29629 \!\cdot\! 36643  $  &  $2^{13}\!\cdot\!      29633  \!\cdot\! 36637  $  &  $2^{13}\!\cdot\!      29683 \!\cdot\! 36563   $\\  
$2^{13}\!\cdot\!29741 \!\cdot\! 36473   $  &  $2^{13}\!\cdot\!     29753 \!\cdot\! 36457  $  &  $2^{13}\!\cdot\!      29803 \!\cdot\!  36383   $  &  $2^{13}\!\cdot\!     29989 \!\cdot\! 36107    $\\   
$2^{13}\!\cdot\!29989\!\cdot\!  36109  $  &  $2^{13}\!\cdot\!      30089\!\cdot\!  35963 $  &  $2^{13}\!\cdot\!      30091\!\cdot\!   35963   $  &  $2^{13}\!\cdot\!     30119\!\cdot\!  35923     $\\ 
$2^{13}\!\cdot\!30137 \!\cdot\! 35897   $  &  $2^{13}\!\cdot\!     30161 \!\cdot\! 35863   $  &  $2^{13}\!\cdot\!     30223 \!\cdot\!  35771    $  &  $2^{13}\!\cdot\!    30253 \!\cdot\! 35729   $\\    
$2^{13}\!\cdot\!30391 \!\cdot\! 35543   $  &  $2^{13}\!\cdot\!     30403 \!\cdot\! 35527  $  &  $2^{13}\!\cdot\!      30469 \!\cdot\!  35437   $  &  $2^{13}\!\cdot\!     30491 \!\cdot\! 35407   $\\    
$2^{13}\!\cdot\!30509 \!\cdot\! 35381    $  &  $2^{13}\!\cdot\!    30697 \!\cdot\! 35129  $  &  $2^{13}\!\cdot\!      30757  \!\cdot\! 35053   $  &  $2^{13}\!\cdot\!     30839  \!\cdot\!34949   $\\
\hline\end{tabular}}

\vfill\eject
\noindent{\bf Table 21.}\quad All primitive weird numbers of form $2^{13}pq$. (Part 9)
\vskip 6pt
\centerline{\begin{tabular}{|llll|}\hline& & & \\  
$2^{13}\!\cdot\!32621 \!\cdot\! 32783  $  &  $2^{13}\!\cdot\!      32621\!\cdot\!  32789  $  &  $2^{13}\!\cdot\!      32621\!\cdot\!   32797  $  &  $2^{13}\!\cdot\!      32621 \!\cdot\! 32801    $\\   
$2^{13}\!\cdot\!32621 \!\cdot\! 32803  $  &  $2^{13}\!\cdot\!      32621 \!\cdot\! 32831 $  &  $2^{13}\!\cdot\!       32621  \!\cdot\! 32833  $  &  $2^{13}\!\cdot\!      32621\!\cdot\!  32839  $\\   
$2^{13}\!\cdot\!32621 \!\cdot\! 32843   $  &  $2^{13}\!\cdot\!     32633  \!\cdot\!32771 $  &  $2^{13}\!\cdot\!       32633 \!\cdot\!  32779    $  &  $2^{13}\!\cdot\!    32633 \!\cdot\! 32783     $\\  
$2^{13}\!\cdot\!32633 \!\cdot\! 32789  $  &  $2^{13}\!\cdot\!      32633  \!\cdot\!32797  $  &  $2^{13}\!\cdot\!      32633 \!\cdot\!  32801   $  &  $2^{13}\!\cdot\!     32633 \!\cdot\! 32803     $\\  
$2^{13}\!\cdot\!32633 \!\cdot\! 32831 $  &  $2^{13}\!\cdot\!       32633 \!\cdot\! 32833   $  &  $2^{13}\!\cdot\!     32633 \!\cdot\!  32839  $  &  $2^{13}\!\cdot\!      32633 \!\cdot\! 32843    $\\  
 $2^{13}\!\cdot\!32633 \!\cdot\! 32869  $  &  $2^{13}\!\cdot\!      32647 \!\cdot\! 32771 $  &  $2^{13}\!\cdot\!       32647 \!\cdot\!  32779  $  &  $2^{13}\!\cdot\!      32647 \!\cdot\! 32783     $\\  
$2^{13}\!\cdot\!32647 \!\cdot\! 32789  $  &  $2^{13}\!\cdot\!      32647  \!\cdot\!32797 $  &  $2^{13}\!\cdot\!       32647 \!\cdot\!  32801   $  &  $2^{13}\!\cdot\!     32647\!\cdot\!  32803  $\\
\hline\end{tabular}}

 \vskip 12pt
\noindent{\bf Table 22.}\quad All primitive weird numbers of form $2^{14}pq$. (Part 1)
\vskip 6pt
\centerline{\begin{tabular}{|llll|}\hline& & & \\  
$2^{14}\!\cdot\!32771\!\cdot\!  215256313$  &   $2^{14}\!\cdot\!  32771\!\cdot\! 215911733$    &    $2^{14}\!\cdot\!  2771\!\cdot\!   215952697$   &   $2^{14}\!\cdot\! 32771 \!\cdot\!  216108359$\\    
$2^{14}\!\cdot\!32771\!\cdot\!  216132937$ &  $2^{14}\!\cdot\!   32771\!\cdot\!  216345949$ &   $2^{14}\!\cdot\!   32771\!\cdot\!  216411491$ &   $2^{14}\!\cdot\!    32771\!\cdot\!  216501611  $\\           
$2^{14}\!\cdot\!32771\!\cdot\!  216894863$ &   $2^{14}\!\cdot\!     32771\!\cdot\!  216919441$ &  $2^{14}\!\cdot\!     32771\!\cdot\!  217091489$ &  $2^{14}\!\cdot\!    32771\!\cdot\!  217132453 $\\            
$2^{14}\!\cdot\!32771\!\cdot\!  217771487$ &  $2^{14}\!\cdot\!     32771\!\cdot\!  217877993$ & $2^{14}\!\cdot\!        32771\!\cdot\!  217902571$ &  $2^{14}\!\cdot\!   32771\!\cdot\!  217984499$\\
$2^{14}\!\cdot\!32771\!\cdot\!  218205703$ &  $2^{14}\!\cdot\!     32771\!\cdot\!  218312209 $ &  $2^{14}\!\cdot\!      32771\!\cdot\!  218402329$ &    $2^{14}\!\cdot\!   32771\!\cdot\!  218754617 $\\
$2^{14}\!\cdot\!32771\!\cdot\!  218795581$ &  $2^{14}\!\cdot\!   32771\!\cdot\!  218902087$ & $2^{14}\!\cdot\!     32771\!\cdot\!  219278953$ &   $2^{14}\!\cdot\!    32771\!\cdot\!  219541121$\\
 $2^{14}\!\cdot\!32771\!\cdot\!  219844253$ &  $2^{14}\!\cdot\!  32771\!\cdot\!  219885217$ &   $2^{14}\!\cdot\!    32771\!\cdot\!  220262083$ &   $2^{14}\!\cdot\!     32771 \!\cdot\! 220458709 $\\
\hline\end{tabular}}

\vfill\eject
\noindent{\bf Table 23.}\quad All primitive weird numbers of form $2^{14}pq$. (Part 2)
\vskip 6pt
\centerline{\begin{tabular}{|llll|}\hline& & & \\     
$2^{14}\!\cdot\!32771\!\cdot\!224563277$ & $2^{14}\!\cdot\!32771\!\cdot\!224866409$ & $2^{14}\!\cdot\!32771\!\cdot\!225284239$ & $2^{14}\!\cdot\!32771\!\cdot\!226177249$\\          
$2^{14}\!\cdot\!32771\!\cdot\!226357489$ & $2^{14}\!\cdot\!32771\!\cdot\!226570501$ & $2^{14}\!\cdot\!32771\!\cdot\!226816283$ & $2^{14}\!\cdot\!32771\!\cdot\!226947367$\\          
$2^{14}\!\cdot\!32771\!\cdot\!227316041$ & $2^{14}\!\cdot\!32771\!\cdot\!227553631$ & $2^{14}\!\cdot\!32771\!\cdot\!227643751$ & $2^{14}\!\cdot\!32771\!\cdot\!227815799$\\
$2^{14}\!\cdot\!32771\!\cdot\!227930497$ & $2^{14}\!\cdot\!32771\!\cdot\!228323749$ & $2^{14}\!\cdot\!32771\!\cdot\!228495797$ & $2^{14}\!\cdot\!32771\!\cdot\!228536761$\\          
$2^{14}\!\cdot\!32771\!\cdot\!228798929$ & $2^{14}\!\cdot\!32771\!\cdot\!228889049$ & $2^{14}\!\cdot\!32771\!\cdot\!228913627$ & $2^{14}\!\cdot\!32771\!\cdot\!229282301$\\          
$2^{14}\!\cdot\!32771\!\cdot\!229372421$ & $2^{14}\!\cdot\!32771\!\cdot\!229700131$ & $2^{14}\!\cdot\!32771\!\cdot\!229765673$ & $2^{14}\!\cdot\!32771\!\cdot\!229896757$\\          
$2^{14}\!\cdot\!32771\!\cdot\!230109769$ & $2^{14}\!\cdot\!32771\!\cdot\!230371937$ & $2^{14}\!\cdot\!32771\!\cdot\!230593141$ & $2^{14}\!\cdot\!32771\!\cdot\!230699647$\\
\hline\end{tabular}}
\vfill\eject
\noindent{\bf Table 24.}\quad All primitive weird numbers of form $2^{14}pq$. (Part 3)
\vskip 6pt
\centerline{\begin{tabular}{|llll|}\hline& & & \\ 
$2^{14}\!\cdot\!32771\!\cdot\!261635471$ & $2^{14}\!\cdot\!32771\!\cdot\!261750169$ & $2^{14}\!\cdot\!32771\!\cdot\!262225349$ & $2^{14}\!\cdot\!32771\!\cdot\!262356433$\\
$2^{14}\!\cdot\!32771\!\cdot\!262405589$ & $2^{14}\!\cdot\!32771\!\cdot\!262553059$ & $2^{14}\!\cdot\!32771\!\cdot\!262946311$ & $2^{14}\!\cdot\!32771\!\cdot\!263101973$\\          
$2^{14}\!\cdot\!32771\!\cdot\!263126551$ & $2^{14}\!\cdot\!32771\!\cdot\!263691851$ & $2^{14}\!\cdot\!32771\!\cdot\!264306307$ & $2^{14}\!\cdot\!32771\!\cdot\!264519319$\\          
$2^{14}\!\cdot\!32771\!\cdot\!264978113$ & $2^{14}\!\cdot\!32771\!\cdot\!265109197$ & $2^{14}\!\cdot\!32771\!\cdot\!265158353$ & $2^{14}\!\cdot\!32771\!\cdot\!265199317$\\          
$2^{14}\!\cdot\!32771\!\cdot\!265502449$ & $2^{14}\!\cdot\!32771\!\cdot\!265682689$ & $2^{14}\!\cdot\!32771\!\cdot\!265944857$ & $2^{14}\!\cdot\!32771\!\cdot\!265985821$\\          
$2^{14}\!\cdot\!32771\!\cdot\!266075941$ & $2^{14}\!\cdot\!32771\!\cdot\!266141483$ & $2^{14}\!\cdot\!32771\!\cdot\!266157869$ & $2^{14}\!\cdot\!32771\!\cdot\!266837867$\\
$2^{14}\!\cdot\!32771\!\cdot\!266968951$ & $2^{14}\!\cdot\!32771\!\cdot\!267059071$ & $2^{14}\!\cdot\!32771\!\cdot\!267362203$ & $2^{14}\!\cdot\!32771\!\cdot\!267648949$\\
\hline\end{tabular}}

\vfill\eject
\noindent{\bf Table 25.}\quad All primitive weird numbers of form $2^{14}pq$. (Part 4)
\vskip 6pt
\centerline{\begin{tabular}{|llll|}\hline& & & \\ 
$2^{14}\!\cdot\!32783\!\cdot\!64055933$ & $2^{14}\!\cdot\!32783\!\cdot\!64144037$ & $2^{14}\!\cdot\!32783\!\cdot\!64152233$ & $2^{14}\!\cdot\!32783\!\cdot\!64236239$\\              
$2^{14}\!\cdot\!32783\!\cdot\!64271071$ & $2^{14}\!\cdot\!32783\!\cdot\!64348931$ & $2^{14}\!\cdot\!32783\!\cdot\!64375567$ & $2^{14}\!\cdot\!32783\!\cdot\!64424741$\\              
$2^{14}\!\cdot\!32783\!\cdot\!64428839$ & $2^{14}\!\cdot\!32783\!\cdot\!64432937$ & $2^{14}\!\cdot\!32783\!\cdot\!64459573$ & $2^{14}\!\cdot\!32783\!\cdot\!64463671$\\              
$2^{14}\!\cdot\!32783\!\cdot\!64568167$ & $2^{14}\!\cdot\!32783\!\cdot\!64625537$ & $2^{14}\!\cdot\!32783\!\cdot\!64646027$ & $2^{14}\!\cdot\!32783\!\cdot\!64676761$\\
$2^{14}\!\cdot\!32783\!\cdot\!64738229$ & $2^{14}\!\cdot\!32783\!\cdot\!64777159$ & $2^{14}\!\cdot\!32783\!\cdot\!64830431$ & $2^{14}\!\cdot\!32783\!\cdot\!64852969$\\              
$2^{14}\!\cdot\!32783\!\cdot\!64922633$ & $2^{14}\!\cdot\!32783\!\cdot\!64945171$ & $2^{14}\!\cdot\!32783\!\cdot\!64953367$ & $2^{14}\!\cdot\!32783\!\cdot\!64965661$\\              
$2^{14}\!\cdot\!32783\!\cdot\!65027129$ & $2^{14}\!\cdot\!32783\!\cdot\!65145967$ & $2^{14}\!\cdot\!32783\!\cdot\!65154163$ & $2^{14}\!\cdot\!32783\!\cdot\!65162359$\\
\hline\end{tabular}}

\vfill\eject
\noindent{\bf Table 26.}\quad All primitive weird numbers of form $2^{14}pq$. (Part 5)
\vskip 6pt
\centerline{\begin{tabular}{|llll|}\hline& & & \\ 
$2^{14}\!\cdot\!32797\!\cdot\!34862117$ & $2^{14}\!\cdot\!32797\!\cdot\!34867583$ & $2^{14}\!\cdot\!32797\!\cdot\!34868677$ & $2^{14}\!\cdot\!32797\!\cdot\!34874143$\\            
$2^{14}\!\cdot\!32797\!\cdot\!34899287$ & $2^{14}\!\cdot\!32797\!\cdot\!34920059$ & $2^{14}\!\cdot\!32797\!\cdot\!34930991$ & $2^{14}\!\cdot\!32797\!\cdot\!34941923$\\
$2^{14}\!\cdot\!32797\!\cdot\!34943017$ & $2^{14}\!\cdot\!32797\!\cdot\!34948483$ & $2^{14}\!\cdot\!32797\!\cdot\!34957229$ & $2^{14}\!\cdot\!32797\!\cdot\!34968161$\\              
$2^{14}\!\cdot\!32797\!\cdot\!34979093$ & $2^{14}\!\cdot\!32797\!\cdot\!34985653$ & $2^{14}\!\cdot\!32797\!\cdot\!35010797$ & $2^{14}\!\cdot\!32797\!\cdot\!35012983$\\              
$2^{14}\!\cdot\!32797\!\cdot\!35021729$ & $2^{14}\!\cdot\!32797\!\cdot\!35028289$ & $2^{14}\!\cdot\!32797\!\cdot\!35044687$ & $2^{14}\!\cdot\!32797\!\cdot\!35116841$\\              
$2^{14}\!\cdot\!32797\!\cdot\!35124493$ & $2^{14}\!\cdot\!32797\!\cdot\!35127773$ & $2^{14}\!\cdot\!32797\!\cdot\!35159477$ & $2^{14}\!\cdot\!32797\!\cdot\!35161663$\\              
$2^{14}\!\cdot\!32797\!\cdot\!35170409$ & $2^{14}\!\cdot\!32797\!\cdot\!35230537$ & $2^{14}\!\cdot\!32797\!\cdot\!35244749$ & $2^{14}\!\cdot\!32797\!\cdot\!35340953$\\
\hline\end{tabular}}

\vfill\eject
\noindent{\bf Table 27.}\quad All primitive weird numbers of form $2^{14}pq$. (Part 6)
\vskip 6pt
\centerline{\begin{tabular}{|llll|}\hline& & & \\  
$2^{14}\!\cdot\!32831\!\cdot\!16758173$ & $2^{14}\!\cdot\!32831\!\cdot\!16759199$ & $2^{14}\!\cdot\!32831\!\cdot\!16763303$ & $2^{14}\!\cdot\!32831\!\cdot\!16767407$\\              
$2^{14}\!\cdot\!32831\!\cdot\!16777153$ & $2^{14}\!\cdot\!32831\!\cdot\!16785361$ & $2^{14}\!\cdot\!32831\!\cdot\!16790491$ & $2^{14}\!\cdot\!32831\!\cdot\!16791517$\\              
$2^{14}\!\cdot\!32833\!\cdot\!16062799$ & $2^{14}\!\cdot\!32833\!\cdot\!16067773$ & $2^{14}\!\cdot\!32833\!\cdot\!16068271$ & $2^{14}\!\cdot\!32833\!\cdot\!16075733$\\              
$2^{14}\!\cdot\!32833\!\cdot\!16084687$ & $2^{14}\!\cdot\!32833\!\cdot\!16107571$ & $2^{14}\!\cdot\!32833\!\cdot\!16131449$ & $2^{14}\!\cdot\!32833\!\cdot\!16142891$\\              
$2^{14}\!\cdot\!32833\!\cdot\!16152841$ & $2^{14}\!\cdot\!32833\!\cdot\!16171247$ & $2^{14}\!\cdot\!32833\!\cdot\!16175227$ & $2^{14}\!\cdot\!32833\!\cdot\!16180699$\\
$2^{14}\!\cdot\!32833\!\cdot\!16182191$ & $2^{14}\!\cdot\!32833\!\cdot\!16182689$ & $2^{14}\!\cdot\!32833\!\cdot\!16204577$ & $2^{14}\!\cdot\!32833\!\cdot\!16210547$\\              
$2^{14}\!\cdot\!32833\!\cdot\!16216019$ & $2^{14}\!\cdot\!32833\!\cdot\!16225969$ & $2^{14}\!\cdot\!32833\!\cdot\!16227461$ & $2^{14}\!\cdot\!32833\!\cdot\!16231441$\\
\hline\end{tabular}}

\vfill\eject
\noindent{\bf Table 28.}\quad All primitive weird numbers of form $2^{14}pq$. (Part 7)
\vskip 6pt
\centerline{\begin{tabular}{|llll|}\hline& & & \\ 
$2^{14}\!\cdot\!32917\!\cdot\!7177661$ & $2^{14}\!\cdot\!32917\!\cdot\!7181611$ & $2^{14}\!\cdot\!32917\!\cdot\!7181831$ & $2^{14}\!\cdot\!32917\!\cdot\!7187317$\\  
$2^{14}\!\cdot\!32933\!\cdot\!6462803$ & $2^{14}\!\cdot\!32933\!\cdot\!6465581$ & $2^{14}\!\cdot\!32933\!\cdot\!7465779$ & $2^{14}\!\cdot\!32933\!\cdot\!6466771$\\                   
$2^{14}\!\cdot\!32933\!\cdot\!6466969$ & $2^{14}\!\cdot\!32933\!\cdot\!6472723$ & $2^{14}\!\cdot\!32933\!\cdot\!6475103$ & $2^{14}\!\cdot\!32933\!\cdot\!6475897$\\
$2^{14}\!\cdot\!32933\!\cdot\!6477683$ & $2^{14}\!\cdot\!32933\!\cdot\!6483833$ & $2^{14}\!\cdot\!32933\!\cdot\!6484031$ & $2^{14}\!\cdot\!32933\!\cdot\!6489983$\\                  
$2^{14}\!\cdot\!32933\!\cdot\!6499307$ & $2^{14}\!\cdot\!32933\!\cdot\!6500101$ & $2^{14}\!\cdot\!32939\!\cdot\!6240791$ & $2^{14}\!\cdot\!32939\!\cdot\!6244621$\\  
$2^{14}\!\cdot\!32939\!\cdot\!6245387$ & $2^{14}\!\cdot\!32939\!\cdot\!6249409$ & $2^{14}\!\cdot\!32939\!\cdot\!6253813$ & $2^{14}\!\cdot\!32939\!\cdot\!6254579$\\  
$2^{14}\!\cdot\!32939\!\cdot\!6254771$ & $2^{14}\!\cdot\!32939\!\cdot\!6258793$ & $2^{14}\!\cdot\!32939\!\cdot\!6259367$ & $2^{14}\!\cdot\!32939\!\cdot\!6263963$\\
\hline\end{tabular}}

\vfill\eject
\noindent{\bf Table 29.}\quad All primitive weird numbers of form $2^{14}pq$. (Part 8)
\vskip 6pt
\centerline{\begin{tabular}{|llll|}\hline& & & \\ 
$2^{14}\!\cdot\!33049\!\cdot\!3836731$ & $2^{14}\!\cdot\!33049\!\cdot\!3837199$ & $2^{14}\!\cdot\!33049\!\cdot\!3839309$ & $2^{14}\!\cdot\!33053\!\cdot\!3777749$\\
$2^{14}\!\cdot\!33053\!\cdot\!3779021$ & $2^{14}\!\cdot\!33053\!\cdot\!3780061$ & $2^{14}\!\cdot\!33053\!\cdot\!3780407$ & $2^{14}\!\cdot\!33053\!\cdot\!3782257$\\                   
$2^{14}\!\cdot\!33053\!\cdot\!3783181$ & $2^{14}\!\cdot\!33053\!\cdot\!3785261$ & $2^{14}\!\cdot\!33053\!\cdot\!3786533$ & $2^{14}\!\cdot\!33071\!\cdot\!3555241$\\  
$2^{14}\!\cdot\!33071\!\cdot\!3556981$ & $2^{14}\!\cdot\!33071\!\cdot\!3558613$ & $2^{14}\!\cdot\!33071\!\cdot\!3560027$ & $2^{14}\!\cdot\!33071\!\cdot\!3560789$\\  
$2^{14}\!\cdot\!33071\!\cdot\!3561659$ & $2^{14}\!\cdot\!33071\!\cdot\!3562529$ & $2^{14}\!\cdot\!33071\!\cdot\!3563291$ & $2^{14}\!\cdot\!33073\!\cdot\!3530921$\\                   
$2^{14}\!\cdot\!33073\!\cdot\!3531569$ & $2^{14}\!\cdot\!33073\!\cdot\!3531677$ & $2^{14}\!\cdot\!33073\!\cdot\!3533731$ & $2^{14}\!\cdot\!33073\!\cdot\!3533839$\\
$2^{14}\!\cdot\!33073\!\cdot\!3534379$ & $2^{14}\!\cdot\!33073\!\cdot\!3538919$ & $2^{14}\!\cdot\!33073\!\cdot\!3539027$ & $2^{14}\!\cdot\!33083\!\cdot\!3423881$\\
\hline\end{tabular}}

\vfill\eject
\noindent{\bf Table 30.}\quad All primitive weird numbers of form $2^{14}pq$. (Part 9)
\vskip 6pt
\centerline{\begin{tabular}{|llll|}\hline& & & \\  
$2^{14}\!\cdot\!33377\!\cdot\!1792673$ & $2^{14}\!\cdot\!33391\!\cdot\!1750673$ & $2^{14}\!\cdot\!33391\!\cdot\!1751047$ & $2^{14}\!\cdot\!33391\!\cdot\!1751689$\\ 
$2^{14}\!\cdot\!33391\!\cdot\!1752599$ & $2^{14}\!\cdot\!33403\!\cdot\!1718573$ & $2^{14}\!\cdot\!33403\!\cdot\!1720307$ & $2^{14}\!\cdot\!33403\!\cdot\!1720517$\\  
$2^{14}\!\cdot\!33403\!\cdot\!1720937$ & $2^{14}\!\cdot\!33409\!\cdot\!1703963$ & $2^{14}\!\cdot\!33409\!\cdot\!1704067$ & $2^{14}\!\cdot\!33409\!\cdot\!1704119$\\                   
$2^{14}\!\cdot\!33409\!\cdot\!1704431$ & $2^{14}\!\cdot\!33409\!\cdot\!1704587$ & $2^{14}\!\cdot\!33409\!\cdot\!1704847$ & $2^{14}\!\cdot\!33413\!\cdot\!1692683$\\
$2^{14}\!\cdot\!33413\!\cdot\!1692839$ & $2^{14}\!\cdot\!33413\!\cdot\!1693201$ & $2^{14}\!\cdot\!33413\!\cdot\!1693511$ & $2^{14}\!\cdot\!33413\!\cdot\!1694701$\\                   
$2^{14}\!\cdot\!33427\!\cdot\!1658029$ & $2^{14}\!\cdot\!33427\!\cdot\!1658941$ & $2^{14}\!\cdot\!33457\!\cdot\!1587413$ & $2^{14}\!\cdot\!33457\!\cdot\!1588043$\\  
$2^{14}\!\cdot\!33457\!\cdot\!1588189$ & $2^{14}\!\cdot\!33457\!\cdot\!1588577$ & $2^{14}\!\cdot\!33457\!\cdot\!1588819$ & $2^{14}\!\cdot\!33461\!\cdot\!1577729$\\
\hline\end{tabular}}

\vfill\eject
\noindent{\bf Table 31.}\quad All primitive weird numbers of form $2^{14}pq$. (Part 10)
\vskip 6pt
\centerline{\begin{tabular}{|llll|}\hline& & & \\  
$2^{14}\!\cdot\!33931\!\cdot\!954469$ & $2^{14}\!\cdot\!33937\!\cdot\!949771$ & $2^{14}\!\cdot\!33937\!\cdot\!950177$ & $2^{14}\!\cdot\!33967\!\cdot\!927497$\\
$2^{14}\!\cdot\!33997\!\cdot\!905011$ & $2^{14}\!\cdot\!33997\!\cdot\!905453$ & $2^{14}\!\cdot\!34019\!\cdot\!889901$ & $2^{14}\!\cdot\!34039\!\cdot\!876263$\\        
$2^{14}\!\cdot\!34039\!\cdot\!876611$ & $2^{14}\!\cdot\!34057\!\cdot\!864757$ & $2^{14}\!\cdot\!34061\!\cdot\!851603$ & $2^{14}\!\cdot\!34061\!\cdot\!862159$\\                        
$2^{14}\!\cdot\!34061\!\cdot\!862343$ & $2^{14}\!\cdot\!34061\!\cdot\!862369$ & $2^{14}\!\cdot\!34061\!\cdot\!862501$ & $2^{14}\!\cdot\!34123\!\cdot\!824437$\\                        
$2^{14}\!\cdot\!34123\!\cdot\!824513$ & $2^{14}\!\cdot\!34123\!\cdot\!824563$ & $2^{14}\!\cdot\!34127\!\cdot\!818897$ & $2^{14}\!\cdot\!34127\!\cdot\!818947$\\        
$2^{14}\!\cdot\!34127\!\cdot\!819073$ & $2^{14}\!\cdot\!34127\!\cdot\!819173$ & $2^{14}\!\cdot\!34127\!\cdot\!821833$ & $2^{14}\!\cdot\!34129\!\cdot\!819271$\\                        
$2^{14}\!\cdot\!34129\!\cdot\!820649$ & $2^{14}\!\cdot\!34157\!\cdot\!804901$ & $2^{14}\!\cdot\!34157\!\cdot\!805073$ & $2^{14}\!\cdot\!34157\!\cdot\!805097$\\
\hline\end{tabular}}

\vfill\eject
\noindent{\bf Table 32.}\quad All primitive weird numbers of form $2^{14}pq$. (Part 11)
\vskip 6pt
\centerline{\begin{tabular}{|llll|}\hline& & & \\  
$2^{14}\!\cdot\!34961\!\cdot\!521929$ & $2^{14}\!\cdot\!34963\!\cdot\!521483$ & $2^{14}\!\cdot\!34981\!\cdot\!517589$ & $2^{14}\!\cdot\!35051\!\cdot\!502699$\\
$2^{14}\!\cdot\!35053\!\cdot\!502441$ & $2^{14}\!\cdot\!35059\!\cdot\!501013$ & $2^{14}\!\cdot\!35059\!\cdot\!501043$ & $2^{14}\!\cdot\!35059\!\cdot\!501089$\\                        
$2^{14}\!\cdot\!35083\!\cdot\!496297$ & $2^{14}\!\cdot\!35099\!\cdot\!491581$ & $2^{14}\!\cdot\!35099\!\cdot\!491611$ & $2^{14}\!\cdot\!35099\!\cdot\!493177$\\        
$2^{14}\!\cdot\!35141\!\cdot\!484987$ & $2^{14}\!\cdot\!35149\!\cdot\!483409$ & $2^{14}\!\cdot\!35153\!\cdot\!482683$ & $2^{14}\!\cdot\!35171\!\cdot\!479299$\\                        
$2^{14}\!\cdot\!35171\!\cdot\!479387$ & $2^{14}\!\cdot\!35201\!\cdot\!473723$ & $2^{14}\!\cdot\!35201\!\cdot\!473839$ & $2^{14}\!\cdot\!35221\!\cdot\!470201$\\        
$2^{14}\!\cdot\!35257\!\cdot\!463849$ & $2^{14}\!\cdot\!35257\!\cdot\!463891$ & $2^{14}\!\cdot\!35257\!\cdot\!463919$ & $2^{14}\!\cdot\!35267\!\cdot\!462181$\\                        
$2^{14}\!\cdot\!35279\!\cdot\!458879$ & $2^{14}\!\cdot\!35279\!\cdot\!458921$ & $2^{14}\!\cdot\!35279\!\cdot\!460073$ & $2^{14}\!\cdot\!35279\!\cdot\!460087$\\
\hline\end{tabular}}

\vfill\eject
\noindent{\bf Table 33.}\quad All primitive weird numbers of form $2^{14}pq$. (Part 12)
\vskip 6pt
\centerline{\begin{tabular}{|llll|}\hline& & & \\  
$2^{14}\!\cdot\!36857\!\cdot\!295081$ & $2^{14}\!\cdot\!36857\!\cdot\!295153$ & $2^{14}\!\cdot\!36877\!\cdot\!293957$ & $2^{14}\!\cdot\!36887\!\cdot\!293269$\\
$2^{14}\!\cdot\!36913\!\cdot\!291689$ & $2^{14}\!\cdot\!36919\!\cdot\!291359$ & $2^{14}\!\cdot\!36923\!\cdot\!291101$ & $2^{14}\!\cdot\!36929\!\cdot\!290711$\\                        
$2^{14}\!\cdot\!36997\!\cdot\!286547$ & $2^{14}\!\cdot\!36997\!\cdot\!286591$ & $2^{14}\!\cdot\!37013\!\cdot\!285599$ & $2^{14}\!\cdot\!37019\!\cdot\!285251$\\        
$2^{14}\!\cdot\!37061\!\cdot\!282773$ & $2^{14}\!\cdot\!37087\!\cdot\!281251$ & $2^{14}\!\cdot\!37097\!\cdot\!280703$ & $2^{14}\!\cdot\!37097\!\cdot\!280729$\\                        
$2^{14}\!\cdot\!37139\!\cdot\!278321$ & $2^{14}\!\cdot\!37139\!\cdot\!278347$ & $2^{14}\!\cdot\!37201\!\cdot\!274871$ & $2^{14}\!\cdot\!37223\!\cdot\!273709$\\                        
$2^{14}\!\cdot\!37273\!\cdot\!271027$ & $2^{14}\!\cdot\!37277\!\cdot\!270791$ & $2^{14}\!\cdot\!37277\!\cdot\!270799$ & $2^{14}\!\cdot\!37307\!\cdot\!269251$\\        
$2^{14}\!\cdot\!37321\!\cdot\!268507$ & $2^{14}\!\cdot\!37339\!\cdot\!267587$ & $2^{14}\!\cdot\!37357\!\cdot\!266683$ & $2^{14}\!\cdot\!37361\!\cdot\!266447$\\
\hline\end{tabular}}

\vfill\eject
\noindent{\bf Table 34.}\quad All primitive weird numbers of form $2^{14}pq$. (Part 13)
\vskip 6pt
\centerline{\begin{tabular}{|llll|}\hline& & & \\  
$2^{14}\!\cdot\!39293\!\cdot\!196717$ & $2^{14}\!\cdot\!39301\!\cdot\!196583$ & $2^{14}\!\cdot\!39301\!\cdot\!196613$ & $2^{14}\!\cdot\!39301\!\cdot\!196643$\\
$2^{14}\!\cdot\!39301\!\cdot\!196661$ & $2^{14}\!\cdot\!39301\!\cdot\!196709$ & $2^{14}\!\cdot\!39301\!\cdot\!196727$ & $2^{14}\!\cdot\!39301\!\cdot\!196739$\\                        
$2^{14}\!\cdot\!39301\!\cdot\!196751$ & $2^{14}\!\cdot\!39301\!\cdot\!196769$ & $2^{14}\!\cdot\!39301\!\cdot\!196799$ & $2^{14}\!\cdot\!39301\!\cdot\!196817$\\                        
$2^{14}\!\cdot\!39301\!\cdot\!196853$ & $2^{14}\!\cdot\!39301\!\cdot\!196871$ & $2^{14}\!\cdot\!39301\!\cdot\!196901$ & $2^{14}\!\cdot\!39301\!\cdot\!196907$\\        
$2^{14}\!\cdot\!39313\!\cdot\!196583$ & $2^{14}\!\cdot\!39313\!\cdot\!196613$ & $2^{14}\!\cdot\!39313\!\cdot\!196643$ & $2^{14}\!\cdot\!39313\!\cdot\!196661$\\                        
$2^{14}\!\cdot\!39313\!\cdot\!196709$ & $2^{14}\!\cdot\!39313\!\cdot\!196727$ & $2^{14}\!\cdot\!39313\!\cdot\!196739$ & $2^{14}\!\cdot\!39313\!\cdot\!196751$\\        
$2^{14}\!\cdot\!39313\!\cdot\!196769$ & $2^{14}\!\cdot\!39317\!\cdot\!196657$ & $2^{14}\!\cdot\!39317\!\cdot\!196663$ & $2^{14}\!\cdot\!39317\!\cdot\!196681$\\
\hline\end{tabular}%
}
\vskip 24pt
\noindent{\bf Table 42.}\quad Large weird numbers of form $2^kpq$. (Part 1)
\vskip 6pt
\centerline{\begin{tabular}{|l|l|l|l|}\hline& & & \\         
\# &  $p$, $q$   &  Parameters &  Proved by\\ 
\hline & & &\\
1.  (219 digits)& $ p=2^{241}+2^{84}\Phi_1 \Phi_2\Phi_4\Phi_8\Phi_{24}\Phi_{24}-1$ (90 digits) & $P=3$, $\Delta=13$ & Thm. 4.2.8 (Prop. 2)\\
 $j=10$, $k=240$       & $ q=2^{241}+2^{158}\Phi_3\Phi_6\Phi_{12}-1$ (73 digits) & $P=1$, $\Delta=5$ & Thm. 4.2.3 (Prop. 1)\\  
& & &\\
\hline & & &\\
2.  (234 digits)& $ p=2^{241}+2^{87}\Phi_1 \Phi_2\Phi_4\Phi_6\Phi_8-1$ (73 digits) & $P=12$, $\Delta=148$ & Thm. 4.2.8 (Prop. 2)\\
\hline \end{tabular}}
\vfill\eject

\noindent{\bf Table 43.}\quad Large weird numbers of form $2^kpq$. (Part 2)
\vskip 6pt
\centerline{\begin{tabular}{|l|l|l|l|}\hline& & & \\         
\# &  $p$, $q$   &  Parameters &  Proved by\\ 
\hline & & &\\
7.  (353 digits)             &           $ p=2^{337}+2^{120}\Phi_3 \Phi_6-1$ (102 digits)                                                                       & $P=5$, $\Delta=29$         & Thm. 4.2.8 (Prop. 2)\\
 $j=14$, $k=336$       &           $ q=2^{337}+2^{218}\Phi_1\Phi_2\Phi_4\Phi_8\Phi_{12}\Phi_{24}-1$ (150 digits)               & $P=3$, $\Delta=13$         & Thm. 4.2.8 (Prop. 2)\\  
& & &\\
\hline & & &\\
8.  (368 digits)             &          $ p=2^{361}+2^{170}\Phi_1 \Phi_2\Phi_3\Phi_4\Phi_8\Phi_{12}\Phi_{24}-1$ (151 digits)                           & $P=13$, $\Delta=173$     & Thm. 4.2.8 (Prop. 2)\\
\hline \end{tabular}}

\vfill\eject

\noindent{\bf Table 44.}\quad Large weird numbers of form $2^kpq$. (Part 3)
\vskip 6pt
\centerline{%
}
\vfill\eject

\noindent{\bf Table 45.}\quad Large weird numbers of form $2^kpq$. (Part 4)
\vskip 6pt
\centerline{%
}

    \end{document}